\documentclass[12pt,a4paper]{article}
\usepackage[width=18cm,height=23cm]{geometry}
\usepackage{amsmath,amssymb,amsfonts,amsthm,enumerate,float}

\usepackage{graphicx}
\usepackage{tikz}
\usetikzlibrary{shapes,backgrounds}
\usepackage{tikz-cd}
\usetikzlibrary{%
  matrix,%
  calc,%
  arrows%
}

\usepackage{mathrsfs,amsbsy,bm,mathtools}
\usepackage{indent first}

\newcommand{\wc}{W_{\sf cantor}}
\newcommand{\wa}{W_{\ca}}
\newcommand{\ca}{\mathcal{A}}

\newcommand{\bz}{\mathbb{Z}}
\newcommand{\br}{\mathbb{R}}

\newcommand{\rw}{\rightarrow}

\newtheorem{Thm}{Theorem}[section]%
\newtheorem{Lem}[Thm]{Lemma}%
\newtheorem{Cor}[Thm]{Corollary}%
\theoremstyle{definition}
\newtheorem{Def}[Thm]{Definition}%

\title{\bf Minkowski dimension of the boundaries of the lakes of Wada}
\author{Zhangchi Chen}
\date{\today}

\begin{document}
\maketitle

\begin{quote}
{\small 
\begin{center}
Abstract
\end{center}
The lakes of Wada are three disjoint simply connected domains in $S^2$ with the counterintuitive property that they all have the same boundary. The common boundary is a indecomposable continuum. In this article we calculated the Minkowski dimension of such boundaries. The lakes constructed in the standard Cantor way has $\ln(6)/\ln(3)\approx 1.6309$-dimensional boundary, while in general, for any number in $[1,2]$ we can construct lakes with such dimensional boundaries.}
\end{quote}

\medskip
{\bf Keywords:}~
Minkowski dimension,
Lakes of Wada.

\section{Introduction} 
A famous counterintuitive example in planar topology, called the lakes of Wada, is that there exists three disjoint simply connected domains $U_1,U_2,U_3$ in $S^2$ such that $\partial U_1=\partial U_2=\partial U_3$. The classical constructions in \cite[P60]{Yoneyama-1917} and in \cite{Hocking-Young-2017} explain the infinite procedure, while \cite[Sect. 8]{Hubbard-Oberste-Vorth-1993} presents it with an interesting scenario of the perils of philanthropy.

To simplify computations later, we introduce the lakes of Wada on the unit square constructed in the standard Cantor way:

Day 1: we dig a simply connected blue lake $B_1$ by digging $[0,\tfrac{2}{3}]\times[\tfrac{1}{3},\tfrac{2}{3}]$ on the island so that every point on the island is at most $\tfrac{\sqrt{2}}{3}$ far away from the blue lake.

Day 2: we dig a simply connected red lake $R_2$ on the island with width $\tfrac{1}{9}$ so that every point on the island is at most $\tfrac{\sqrt{2}}{9}$ far away from the red lake.

Day 3: we dig a simply connected green lake $G_3$ on the island with width $\tfrac{1}{27}$ so that every point on the island is at most $\tfrac{\sqrt{2}}{27}$ far away from the green lake.

\begin{minipage}{\linewidth}
      \centering
      \begin{minipage}{0.3\linewidth}
          \begin{figure}[H]
              \includegraphics[width=\linewidth]{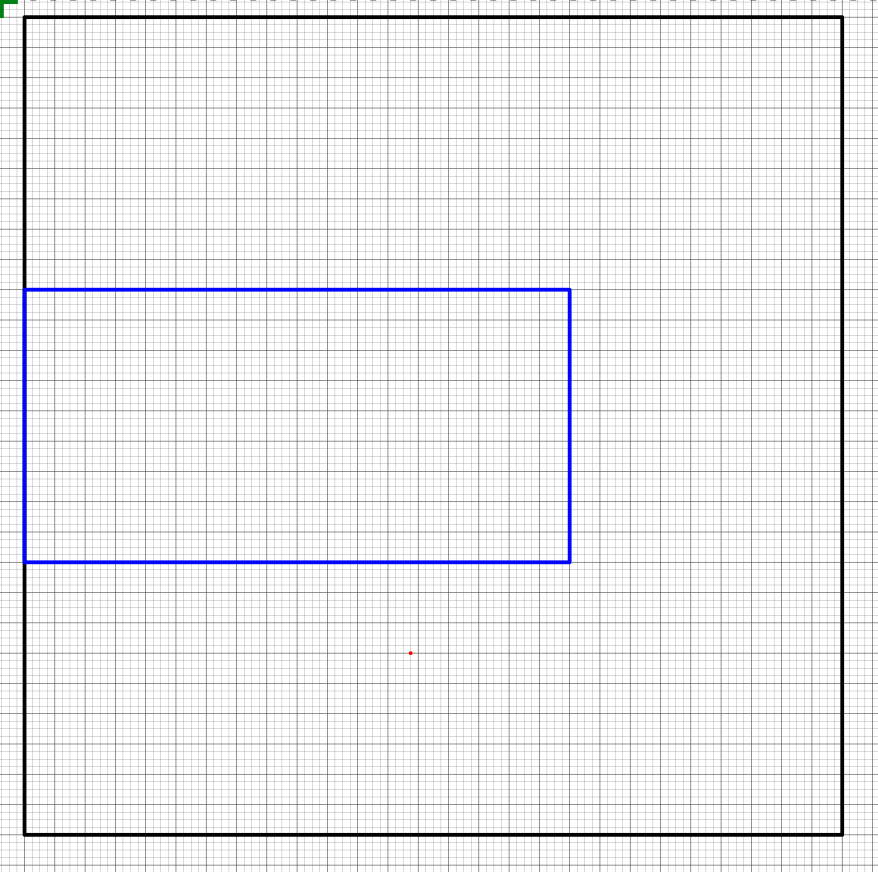}
              \caption{Day 1}
          \end{figure}
      \end{minipage}
      \begin{minipage}{0.3\linewidth}
          \begin{figure}[H]
              \includegraphics[width=\linewidth]{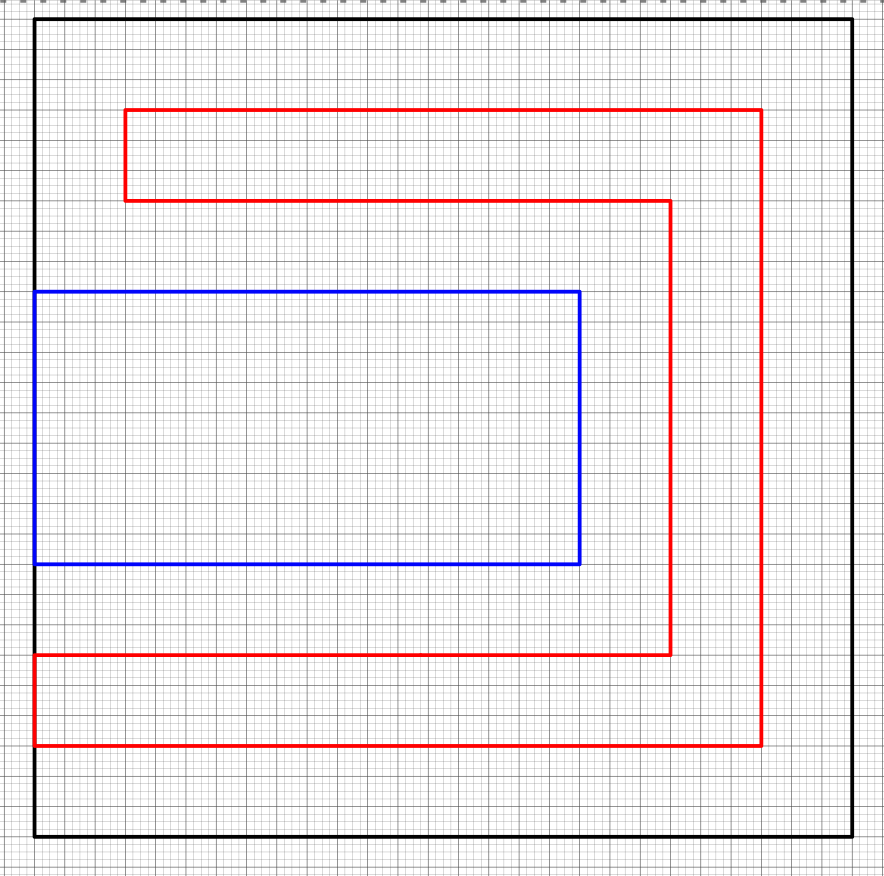}
              \caption{Day 2}
          \end{figure}
      \end{minipage}
      \begin{minipage}{0.3\linewidth}
          \begin{figure}[H]
              \includegraphics[width=\linewidth]{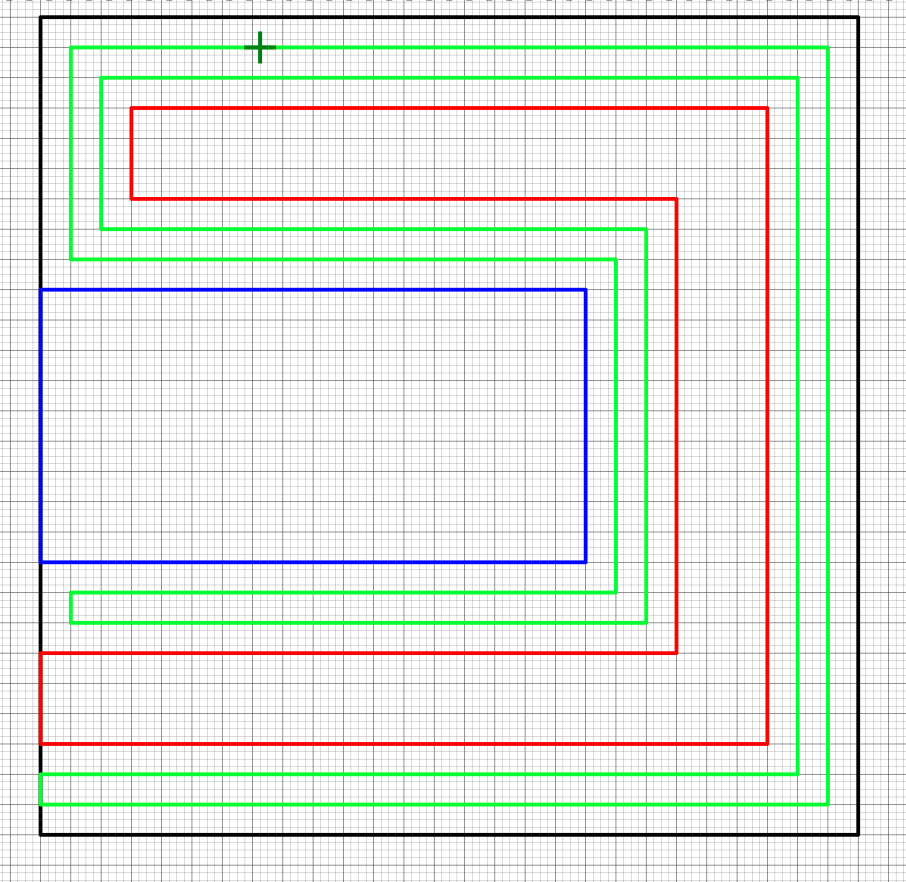}
              \caption{Day 3}
          \end{figure}
      \end{minipage}
\end{minipage}

\medskip

Day 4: we extend our blue lake as $B_4$ by digging a continuous canal of width $\tfrac{1}{81}$ on the island so that every point on the island is at most $\tfrac{\sqrt{2}}{81}$ far away from the blue lake. The canal starts from the peak point of $B_1$: $(\tfrac{2}{3},\tfrac{1}{2})$, separate into two directions: up and down, and terminates at two peak points in violet colour in the lower left corner of the figure.

Day 5: we extend our rad lake as $R_5$ by digging a continuous canal of width $\tfrac{1}{243}$ on the island so that every point on the island is at most $\tfrac{\sqrt{2}}{243}$ far away from the red lake. The canal starts from the peak point of $R_2$, in dark red colour in the upper left corner of the figure, separate into two directions: up and down, and terminates at two peak points, in violet colour near $(\tfrac{2}{3},\tfrac{1}{2})$.

\begin{minipage}{\linewidth}
      \centering
      \begin{minipage}{0.45\linewidth}
          \begin{figure}[H]
              \includegraphics[width=\linewidth]{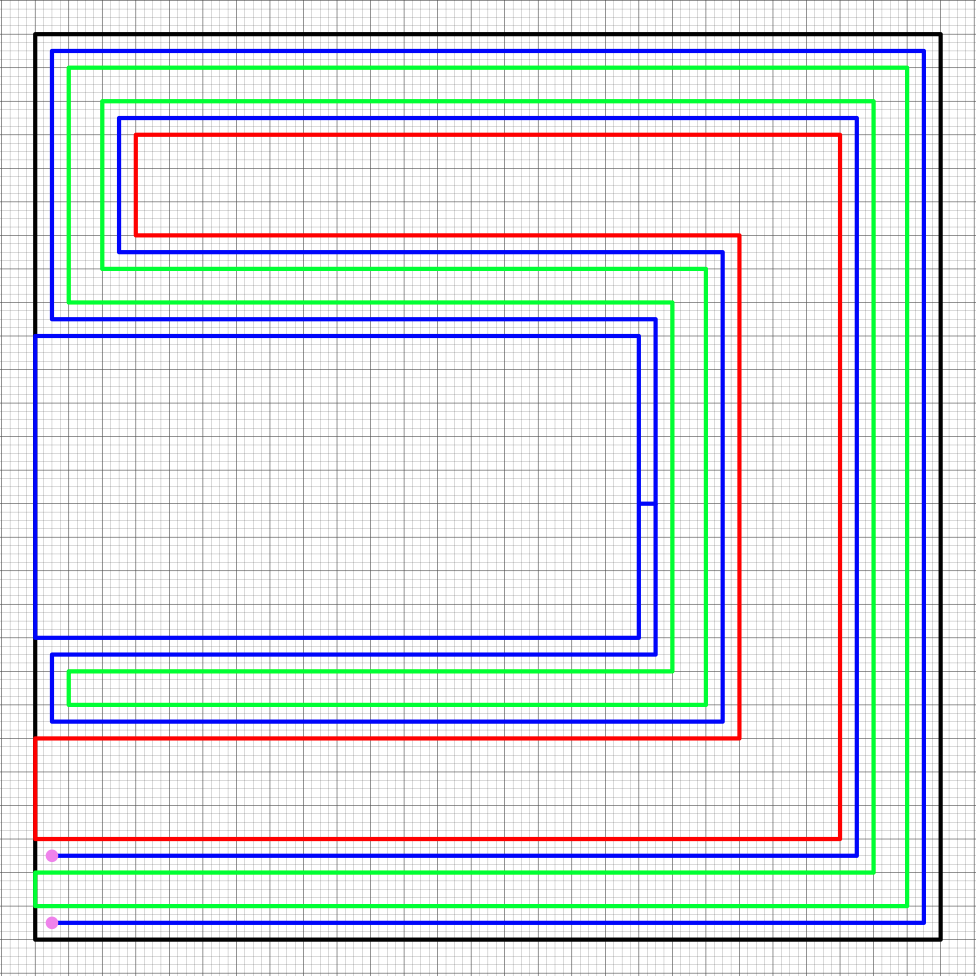}
              \caption{Day 4}
          \end{figure}
      \end{minipage}
      \begin{minipage}{0.45\linewidth}
          \begin{figure}[H]
              \includegraphics[width=\linewidth]{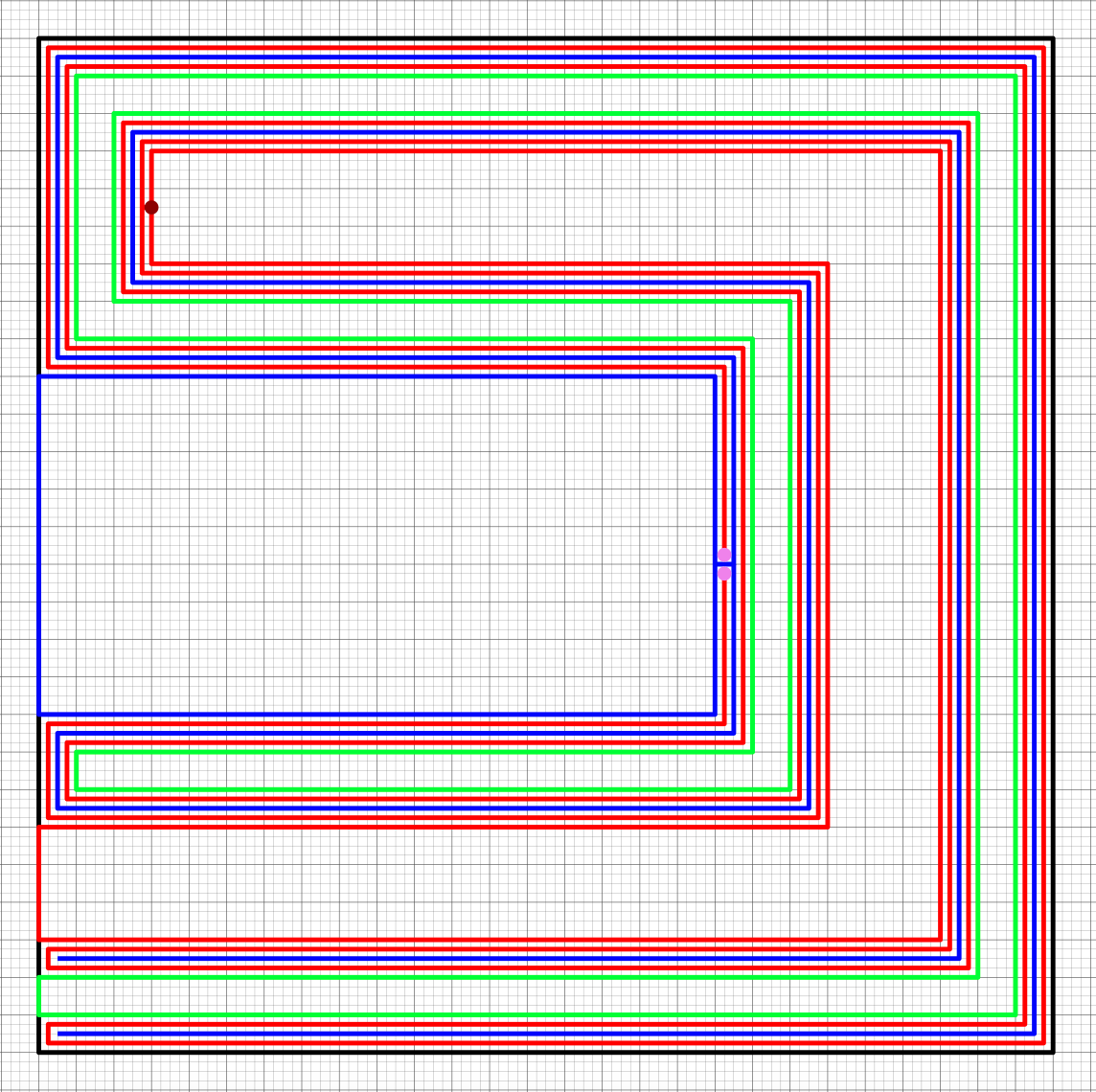}
              \caption{Day 5}
          \end{figure}
      \end{minipage}
\end{minipage}

\medskip

On the following days, we extend the lake with colour green, red, blue in turns, by digging a continuous canal of width $1/3$ of the day before and make sure that every point on the island is closer to our lakes. The canal starts from one of the peak points (at most two), separate into two directions and terminates at two peak points.

Thus we get
\[
\aligned
B_1&\subset B_4\subset B_7\subset\dots\\
R_2&\subset R_5\subset R_8\subset\dots\\
G_3&\subset G_6\subset G_9\subset\dots
\endaligned
\]
three increasing series of domains. The union 
\[
U_B:=\bigcup_{k=0}^\infty B_{3k+1}, \ \ \ U_R:=\bigcup_{k=0}^\infty R_{3k+2}, \ \ \ U_G:=\bigcup_{k=0}^\infty G_{3k+3}
\]
are simply connected domains. Moreover their boundaries in the unit square $(0,1)^2$ are the same, noted by $W_{\sf cantor}:=(0,1)^2\cap\partial U_B=(0,1)^2\cap\partial U_R=(0,1)^2\cap\partial U_G$.

To obtain simply connected domaines on $S^2$ with the same boundary, one can paste the unit square with its reflection, and deform the two squares into two half-balls covering $S^2$. The common edge $V_AV_B$ and $V_A'V_B'$ cover one meridian, while other edges $V_AV_DV_C$ and $V_A'V_D'V_C'$ cover the opposite meridian.

\begin{figure}[H]
\begin{center}
              \includegraphics[width=0.6\linewidth]{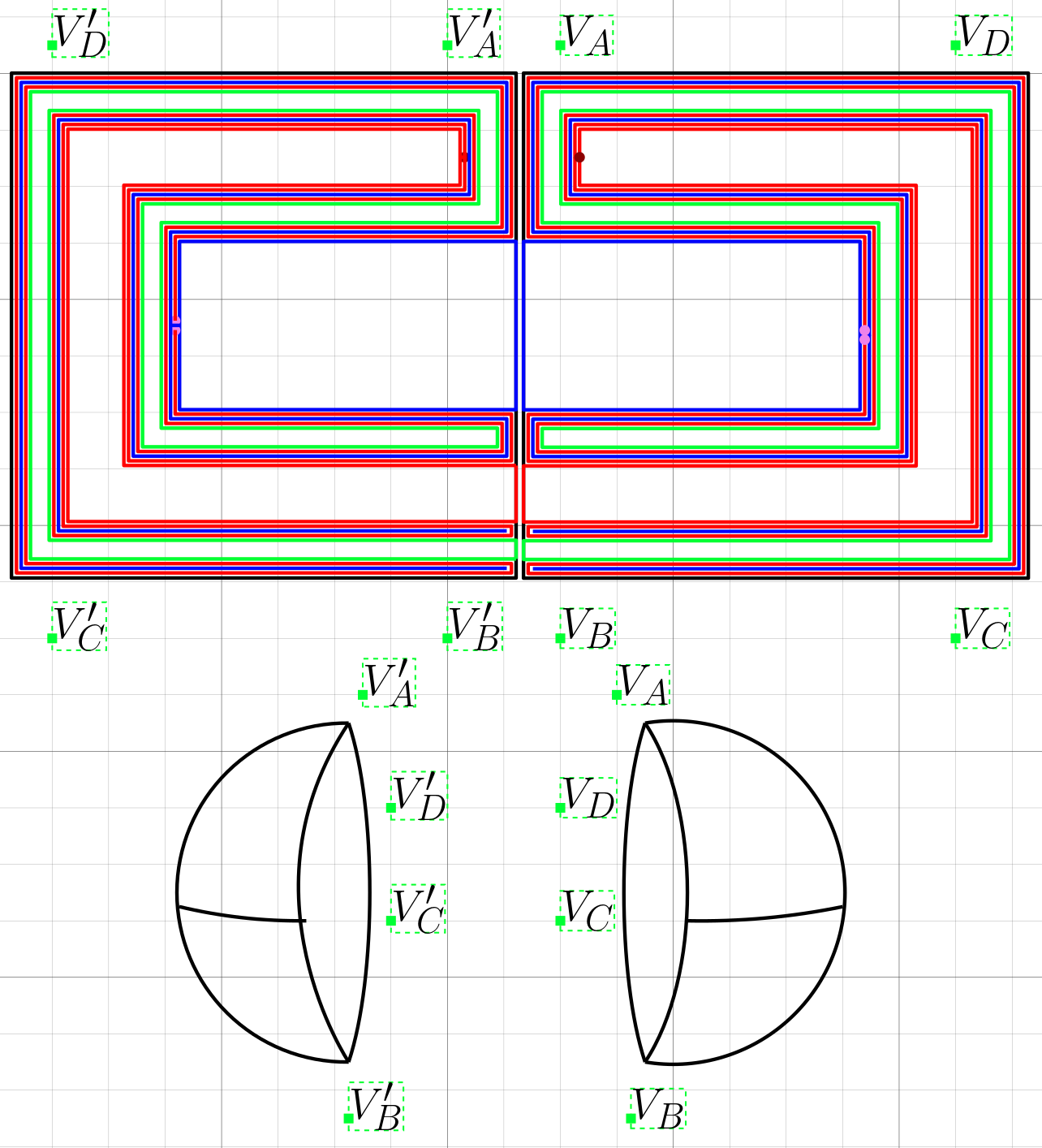}
              \end{center}
              \caption{Gluing two squares on $S^2$}
\end{figure}

The common boundary $W_{\sf cantor}$ is a nowhere dense closed continuous fractal set of infinite length and of 0 area. The Minkowski-Bouligand (also called `box-counting') dimension and Hausdorff-Besicovitch dimension describe how `rough' this fractal set is. We abbreviate them as Minkowski dimension and Hausdorff dimension.

\begin{Def}[Minkowski dimension]
For any bounded set $S\subset \br^n$ and for any $\epsilon>0$, let $N(\epsilon)$ be the number of boxes (squares) required to cover $S$. We define the {Minkowski dimension of $S$}
\[
\mathcal{H}^d(S):=\lim\limits_{r\to 0+}\inf\big\{\sum\limits_{i}r_i^d:\text{there is a countable cover of $S$ by balls with radius $r_i\in(0,r)$}\big\},
\]
and we define the {\em upper Minkowski dimension of $S$} as
\[
\overline{\dim_{\sf box}}(S):=\limsup\limits_{\epsilon\to 0+}\frac{
\ln\big(N(\epsilon)\big)}{\ln(1/\epsilon)}
\]
and the {\em lower Minkowski dimension of $S$} as
\[
\underline{\dim_{\sf box}}(S):=\liminf\limits_{\epsilon\to 0+}\frac{
\ln\big(N(\epsilon)\big)}{\ln(1/\epsilon)}.
\]
When the two numbers agree, we define the {\em Minkowski dimension of $S$} as
\[
\dim_{\sf box}(S):=\lim\limits_{\epsilon\to 0+}\frac{
\ln\big(N(\epsilon)\big)}{\ln(1/\epsilon)}.
\]
\end{Def}

\begin{Def}[Hausdorff dimension]
For any bounded set $S\subset \br^n$ and for any $d>0$ we define the {\em $d$-dimensional Hausdorff outer measure of $S$}
\[
\mathcal{H}^d(S):=\lim\limits_{r\to 0+}\inf\big\{\sum\limits_{i}r_i^d:\text{there is a countable cover of $S$ by balls with radius $r_i\in(0,r)$}\big\},
\]
and we define the {\em Hausdorff dimension of $S$}
\[
\dim_H(S):=\inf\{d\geqslant 0~|~\mathcal{H}^d(S)=0\}.
\]
\end{Def}

The Minkowski dimension is relative easy to calculate, while the Hausdorff dimension deals with countable unions better. There are several basic properties of these dimensions. Assuming $S,S'$ have Minkowski dimensions.

\begin{enumerate}
\item If $S$ is a dimension $d\in\mathbb{Z}$ manifold, then $\dim_H(S)=\dim_{\sf box}(S)=d$;
\item If $S'\subset S$ then
$\dim_H(S')\leqslant \dim_H(S)$ and $\dim_{\sf box}(S')\leqslant\dim_{\sf box}(S)$;
\item If $S=\bigcup\limits_{i=1}^n S_i$ then
$\dim_{\sf box}(S)=\max\dim_{\sf box}(S_i)$;
\item If $S=\bigcup\limits_{i=1}^\infty S_i$ then
$\dim_{H}(S)=\sup\dim_{H}(S_i)$;
\item In general $\dim_{H}(S)\leqslant \dim_{\sf box}(S)$ \cite{Bishop-Peres-2017}. The equality holds when $S$ is strictly self-similar.
\end{enumerate}

From the properties we have $1\leqslant\dim_{H}(W_{\sf cantor})\leqslant \dim_{\sf box}(W_{\sf cantor})\leqslant 2$. In Section \ref{sect-standard-cantor} we prove
\begin{Thm} The boundary of the lakes of Wada construced in the standard Cantor way $W_{\sf cantor}$ has
\[
\dim_{\sf box}(W_{\sf cantor})=\frac{\ln(6)}{\ln(3)}\approx 1.6309.
\]
\end{Thm}

By changing parameters we can construct the lakes of Wada in the modified Cantor way. Take a series of real numbers $\mathcal{A}:=\{a_i>2\}_{i=1}^{\infty}$.

Day 1: we dig a simply connected red lake $B_1$ by digging a continuous canal of width $w_1:=1-\frac{2}{a_1}$ from the left on the island so that every point on the island is at most $\tfrac{\sqrt{2}}{a_1}$ far away from the blue lake.

Day 2: we dig a simply connected red lake $R_2$ on the island with width $w_2:=\frac{1}{a_1}(1-\frac{2}{a_2})$ so that every point on the island is at most $\tfrac{\sqrt{2}}{a_1a_2}$ far away from the red lake.

Day 3: we dig a simply connected green lake $G_3$ on the island with width $w_3:=\tfrac{1}{a_1a_2}(1-\frac{2}{a_3})$ so that every point on the island is at most $\tfrac{\sqrt{2}}{a_1a_2a_3}$ far away from the green lake.

In the future days, write
\[
w_n:=\prod\limits_{i=1}^{n-1}a_i^{-1}(1-\tfrac{2}{a_n}), \ \ \ t_n:=\prod\limits_{i=1}^n a_i^{-1}, \ \ \ (n\geqslant 2)
\]
for the width of canals and $\frac{\sqrt{2}}{2}\times$(the maximal distances between points on the islands and the extended canals on Day $n$).

Then for $k=1,2,\dots$, on day $3k+1$, we extend our blue lake as $B_{3k+1}$ by digging a continuous canal of width $w_{3k+1}$ on the island so that every point on the island is at most $\sqrt{2}\,f_{3k+1}$ far away from the blue lake. On day $3k+2$ we extend the red while keeping it simply connected, while on day $3k+3$ we work on the green. Finally we get the boundary of the lakes of Wada $W_{\mathcal{A}}$ depends on the series $\ca$.

We write $W_{c}$ for $W_{\ca}$ when $\ca$ is a constant sequence $c,c,c,\dots$ with $c>2$. In particular $W_{3}=W_{\sf cantor}$. In Section \ref{sect-modify-cantor} we prove

\begin{Thm}\label{thm-wada-arbitary} Using the notations above
\[
\dim_{\sf box}(W_{c})=\frac{\ln(2c)}{\ln(c)}\rw 1, \ \ \ \ (c\rw +\infty)
\]
In general
\[
\dim_{\sf box}(W_{\ca})=1+\lim\limits_{n\rw+\infty}\frac{n\ln(2)}{\sum\limits_{i=1}^n\ln(a_i)}
\]
\end{Thm}

\begin{Cor}\label{cor-arb} For any real number $d\in[1,2]$ there exists a series $\ca:=\{a_i>2\}_{i=1}^{\infty}$ such that
\[
\dim_{\sf box}(W_{\ca})=d.
\]
\end{Cor}

Higher dimensional lakes of Wada can be obtained by taking products of 2-dimensional lakes of Wada with $\br^{n-2}$.
\begin{Cor} In $\br^n$ with $n\geqslant 2$, for any real number $d\in[n-1,n]$ there exists three disjoint domains $U_1,U_2,U_3$ sharing the same boundary
\[
W=\partial U_1=\partial U_2=\partial U_3
\]
and
\[
\dim_{\sf box}(W)=d.
\]
\end{Cor}

The Hausdorff dimension will be calculated in the next version of the paper.

{\bf Acknowledgement:} The author thanks the organisers of meeting GDR (Groupe de Recherche) GDM (Géométrie Différentielle et Mécanique) at La Rochelle 07-09 July 2021, especially Boris Kolev for introducing the lakes of Wada during a chat in a bar after the meeting.

\section{Minkowski dimension of the boundary in the standard Cantor case}\label{sect-standard-cantor}
From the construction
\[
\wc=\bigcup\limits_{k=0}^{\infty}\partial B_{3k+1}=\bigcup\limits_{k=0}^{\infty}\partial R_{3k+2}=\bigcup\limits_{k=0}^{\infty}\partial G_{3k+3}.
\]
To simplify notations, we write $W_{3k+1},$ $W_{3k+2}$ and $W_{3k+3}$ for $\partial B_{3k+1}$, $\partial R_{3k+2}$ and $\partial G_{3k+3}$ respectively. Thus
\[
\wc=\bigcup\limits_{n=0}^{\infty}W_n.
\]

First we calculate its dimension by box-counting.

Day 1: we cover the blue boundary $\partial B_1$ by exactly 1 square of length 1. So $N(1)=1$. This square also covers $\wc$.

Day 2: we cover the red boundary $\partial R_2$ by 7 squares of length $3^{-1}$. So $N(3^{-1})=7$. These squares also covers $\wc$.

Day 3: we cover the green boundary $\partial G_3$ by 43 squares of length $3^{-2}$. So $N(3^{-2})=43$. These squares also covers $\wc$.

\begin{figure}[H]
\begin{center}
              \includegraphics[width=0.6\linewidth]{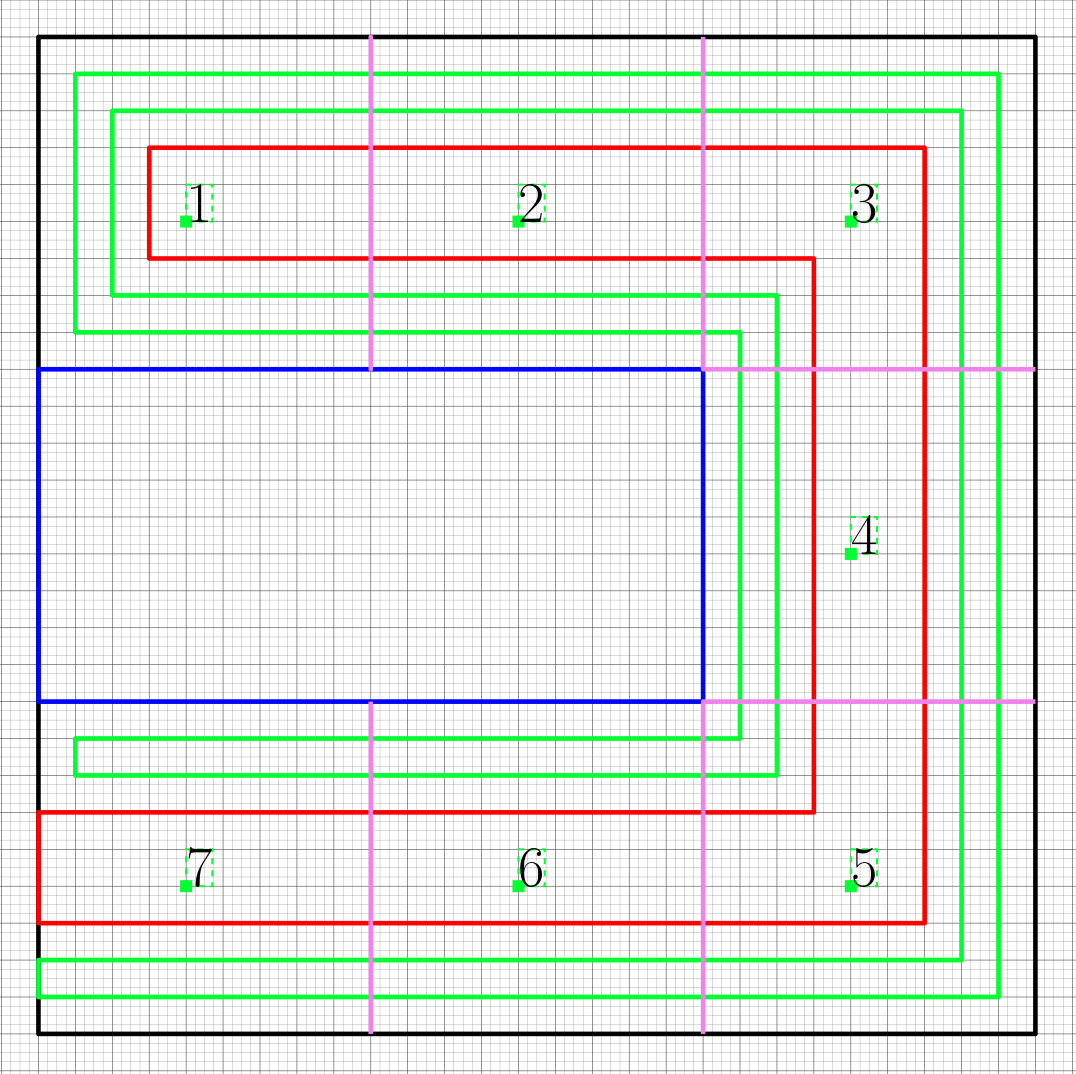}
              \end{center}
              \caption{Day 2: 7 squares of length $3^{-1}$}
\end{figure}

On day $n$, all the $N(3^{1-n})$ squares covering $W_n$ hence covering $\wc$, up to a rotation, are classified into 4 cases:

\begin{itemize}
\item Case 1 (Terminal): The new canal ends in this square, like the square 1 on Day 2.
\item Case 2 (Straight): The new canal passes this square straightly, like squares 2, 4, 6, 7 on Day 2.
\item Case 3 (Turning): The new canal passes this square with a left or right turn, like squares 3, 5 on Day 2.
\item Case 4 (Separation): The new canal enters this square from the left and separates into two branches up and down. It appears exactly once everyday from Day 4.
\end{itemize}

\begin{figure}[H]
\begin{center}
              \includegraphics[width=0.6\linewidth]{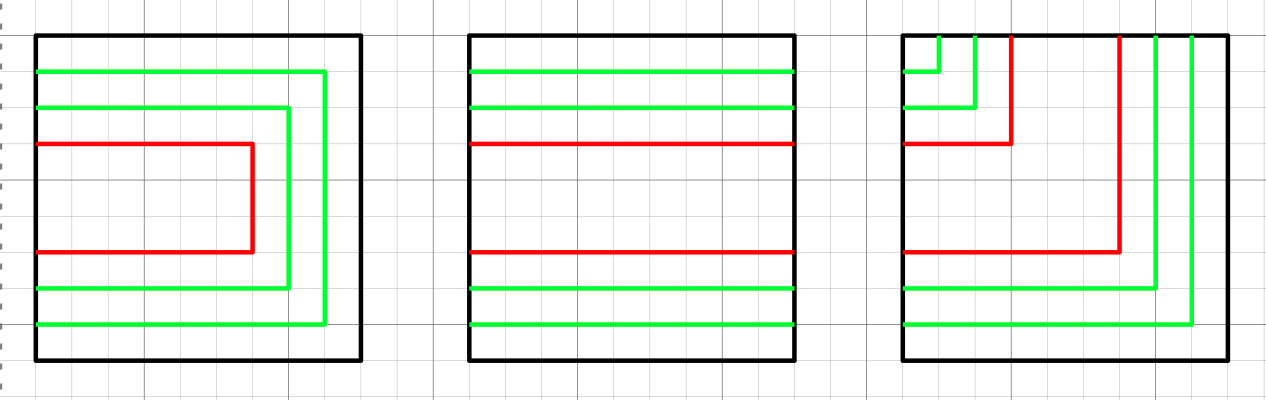}
              \end{center}
              \caption{Terminal, Straight and Turning}
\end{figure}

\begin{figure}[H]
\begin{center}
              \includegraphics[width=0.4\linewidth]{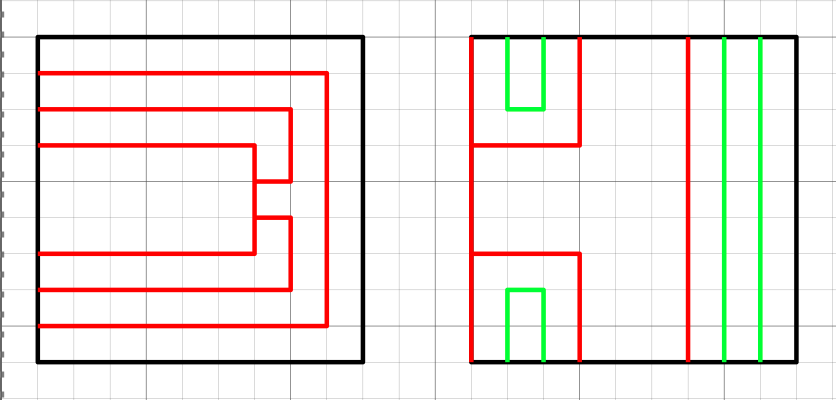}
              \end{center}
              \caption{Separation comes from extending a Terminal, and will split into two Terminals}
\end{figure}

In the first three days, there is exactly one Terminal and no Separation. From Day 4, there will be exactly two Terminals and one Separation.

On day $n+1$, all the $N(3^{-n})$ squares required to cover $W_{n+1}$ are sub-squares of the Terminal, the Straights and the Turnings used one day before. In fact
\begin{itemize}
\item each Terminal is covered by and requires 7 sub-squares;
\item each Straight and each Turning is covered by and requires 6 sub-squares. 
\item each Separation is covered by and requires 5 sub-squares. 
\end{itemize}
thus 
\[
\left\{
\begin{aligned}
N(1)&=1,\\
N(3^{-n})&=6N(3^{1-n})+1, \ \ \ n=1,2,3,\dots
\end{aligned}
\right.
\]
We concludes that $N(3^{-n})=O(6^n)$ so
\[
\lim\limits_{n\rw+\infty}\frac{\ln\big(N(3^{-n})\big)}{\ln(3^n)}=\frac{\ln(6)}{\ln(3)}.
\]
For any $\epsilon\in[3^{-n-1},3^{-n}]$ clearly we have
\[
N(3^{-n})\leqslant N(\epsilon)\leqslant N(3^{-n-1}).
\]
Thus
\[
\aligned
\frac{\ln\big(N(3^{-n})\big)}{\ln(3^{n+1})}&\leqslant \frac{\ln\big(N(\epsilon)\big)}{\ln(1/\epsilon)}\leqslant \frac{\ln\big(N(3^{-n-1})\big)}{\ln(3^{n})}\\
\frac{\ln\big(N(3^{-n})\big)/\ln(3^n)}{1+\ln(3)/\ln(3^n)}&\leqslant \frac{\ln\big(N(\epsilon)\big)}{\ln(1/\epsilon)}\leqslant \frac{\ln\big(N(3^{-n-1})\big)/\ln(3^{n+1})}{1-\ln(3)/\ln(3^{n+1})}
\endaligned
\]
We concludes that the Minkowski dimension exists and
\[
\dim_{\sf box}(\wc)=\frac{\ln(6)}{\ln(3)}\approx 1.6309.
\]

\section{Minkowski dimension of the boundary in the modified Cantor case}\label{sect-modify-cantor}
First we treat the case $\ca:=\{a_i\in\bz_{\geqslant 3}\}_{i=1}^{\infty}$ of integer series. 
\begin{figure}[H]
\begin{center}
              \includegraphics[width=0.6\linewidth]{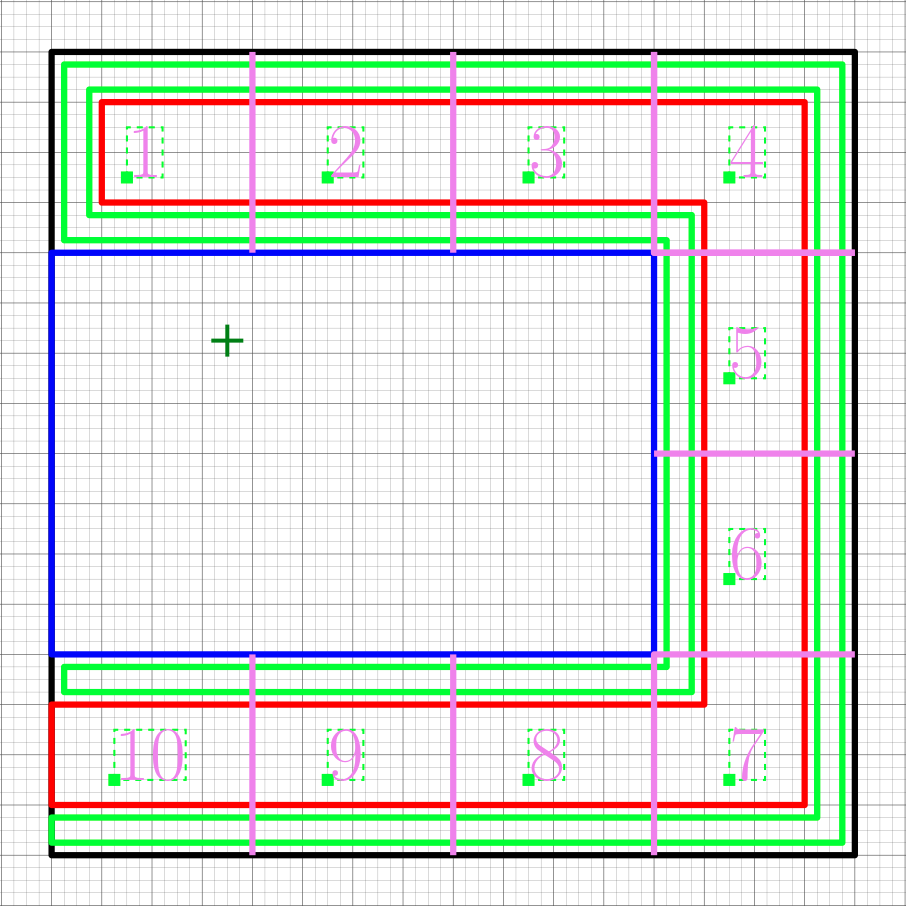}
              \end{center}
              \caption{$W_4$, the lakes of Wada constructed by a constante sequence $4,4,4,\dots$}
\end{figure}

Recall the notations
\[
w_n=\prod\limits_{i=1}^{n-1}a_i^{-1}(1-\tfrac{2}{a_n}), \ \ \ t_n=\prod\limits_{i=1}^n a_i^{-1}, \ \ \ (n\geqslant 2)
\]
for the width of canals and $\frac{\sqrt{2}}{2}\times$(the maximal distances between points on the islands and the extended canals on Day $n$).

On Day $n$ we use $N(t_n)$ squares of length $t_n=\prod\limits_{i=1}^n a_i^{-1}$. They are also classified into one Terminal (two Terminals and one Separation from Day 4) and finitely many Straights and Turnings.

On Day $n+1$ those squares splits into sub-squares of length $t_{n+1}=\prod\limits_{i=1}^{n+1} a_i^{-1}$, while
\begin{itemize}
\item each Terminal is covered by and requires $3a_{n+1}-2$ sub-squares;
\item each Straight and each Turning is covered by and requires $2a_{n+1}$ sub-squares;
\item each Separation is covered by and requires $a_{n+1}+2$ sub-squares.
\end{itemize}
Thus we get recurrence formulas
\[
\left\{
\begin{aligned}
N(t_0)&=N(1)=1,\\
N(t_{n+1})&=2a_{n+1}N({t_n})+a_{n+1}-2, \ \ \ n=1,2,3,\dots
\end{aligned}
\right.
\]
We conclude that $N(t_{n+1})=O\big(\prod\limits_{i=1}^n2a_i\big)$. 

If the sequence $\mathcal{A}$ satisfies
\begin{align}\label{condition-bdd}
\lim\limits_{n\rw+\infty}\frac{\ln(a_{n+1})}{\sum\limits_{i=1}^n\ln(a_i)}=0
\end{align}
for example, if $\mathcal{A}$ is bounded or $a_n=O(n^k)$ for some $k\in\bz$, then by the same argument as before, we conclude that
\[
\dim_{\sf box}(\wa)=\lim\limits_{n\rw+\infty}\frac{\ln\big(N(t_{n})\big)}{\ln(t_n)}=\lim\limits_{n\rw+\infty}\frac{\ln\big(\prod\limits_{i=1}^n2a_i\big)}{\ln\big(\prod\limits_{i=1}^na_i\big)}=1+\lim\limits_{n\rw+\infty}\frac{n\ln(2)}{\sum\limits_{i=1}^na_i}.
\]
But if $\mathcal{A}$ does not satisfy the condition (\ref{condition-bdd}), for example if $a_n=2^{2^n}$, then we need to analyse more delicately.

For any $\epsilon\in(0,1)$ there is a unique $n$ such that $t_{n+1}\leqslant \epsilon<t_n$, i.e.
\[
\prod\limits_{i=1}^n a_i<1/\epsilon\leqslant a_{n+1}\prod\limits_{i=1}^n a_i.
\]
Recall that $a_{n+1}\in\bz_{\geqslant 3}$. Thus there is a unique integer $p\in\{1,2,\dots,a_{n+1}-1\}$ such that
\[
p\prod\limits_{i=1}^n a_i<1/\epsilon\leqslant (p+1)\prod\limits_{i=1}^n a_i.
\]
On Day $n+1$, for any $q\in\{1,2,\dots,a_{n+1}\}$, if we use sub-squares of length $q^{-1}t_n$ instead of $a_{n+1}^{-1}t_n$, then
\begin{itemize}
\item each Terminal is covered by and requires $3q-2$ sub-squares;
\item each Straight and each Turning is covered by and requires $2q$ sub-squares;
\item each Separation is covered by and requires $q+2$ sub-squares.
\end{itemize}
Thus
\[
N(q^{-1}\,t_{n+1})=2qN({t_n})+q-2, \ \ \ q=1,2,\dots,a_{n+1}.
\]
In particular
\[
\aligned
N(p^{-1}t_n)&\leqslant N(\epsilon)\leqslant N\big((p+1)^{-1}t_n\big)\\
2pN(t_n)&\leqslant N(\epsilon)\leqslant 2(p+1)N\big(t_n\big)+p-2<(2p+3)N\big(t_n\big)\\
\frac{\ln(2p)+\ln\big(N(t_n)\big)}{\ln(p+1)+\ln\big(\prod\limits_{i=1}^na_i\big)}&\leqslant \frac{\ln\big(N(\epsilon)\big)}{\ln(\epsilon)}\leqslant \frac{\ln(2p+3)+\ln\big(N(t_n)\big)}{\ln(p)+\ln\big(\prod\limits_{i=1}^na_i\big)}
\endaligned
\]
Both sides converges to $1+\lim\limits_{n\rw+\infty}\frac{n\ln(2)}{\ln(\prod\limits_{i=1}^na_i)}$ The Minkowski dimension is then
\[
\dim_{\sf box}(W_{\ca})=
1+\lim\limits_{n\rw+\infty}\frac{n\ln(2)}{\ln(\prod\limits_{i=1}^na_i)}
\]
In particular, when $\ca$ is a constant sequence $c,c,c,\dots$ we get
\[
\dim_{\sf box}(W_c)=1+\frac{\ln(2)}{\ln(c)}
\]
which tends to 1 when $c\rw+\infty$.


In the general case when $a_n>2$ may not be integers, we may not be able to decompose squares of length $t_n=\prod\limits_{i=1}^{n}a_i^{-1}$ into squares of length $t_{n+1}=\prod\limits_{i=1}^{n+1}a_i^{-1}$. Instead, we can keep square Terminal and Turnings, and rectangular Straights, like the part 1, 2, 3 in the following figures.

\begin{minipage}{\linewidth}
      \centering
      \begin{minipage}{0.45\linewidth}
\begin{figure}[H]
              \includegraphics[width=0.6\linewidth]{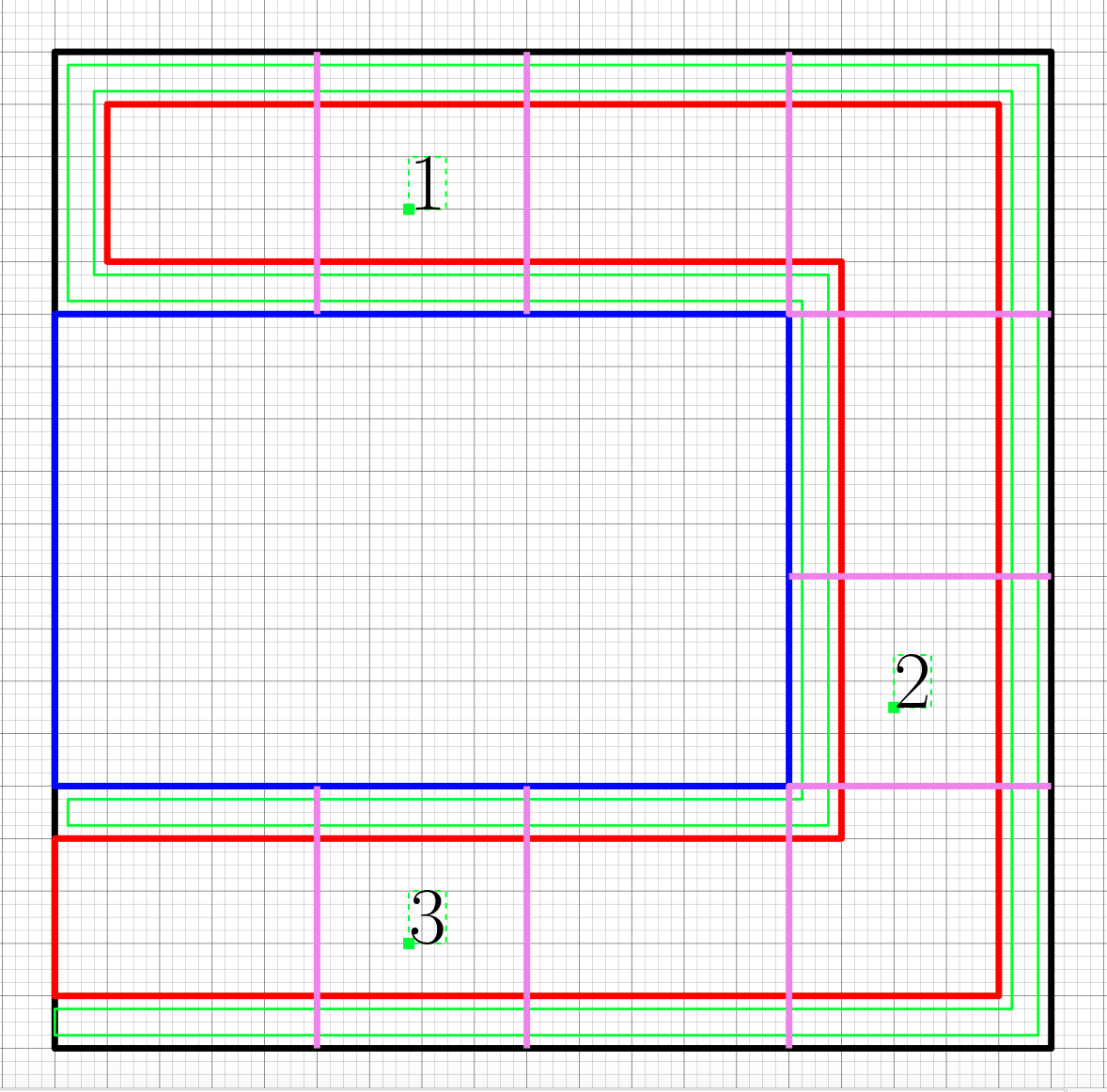}
              \caption{$\wa$ with $\mathcal{A}=3.8,4,4\dots$}
\end{figure}
      \end{minipage}
      \begin{minipage}{0.45\linewidth}
\begin{figure}[H]
              \includegraphics[width=0.6\linewidth]{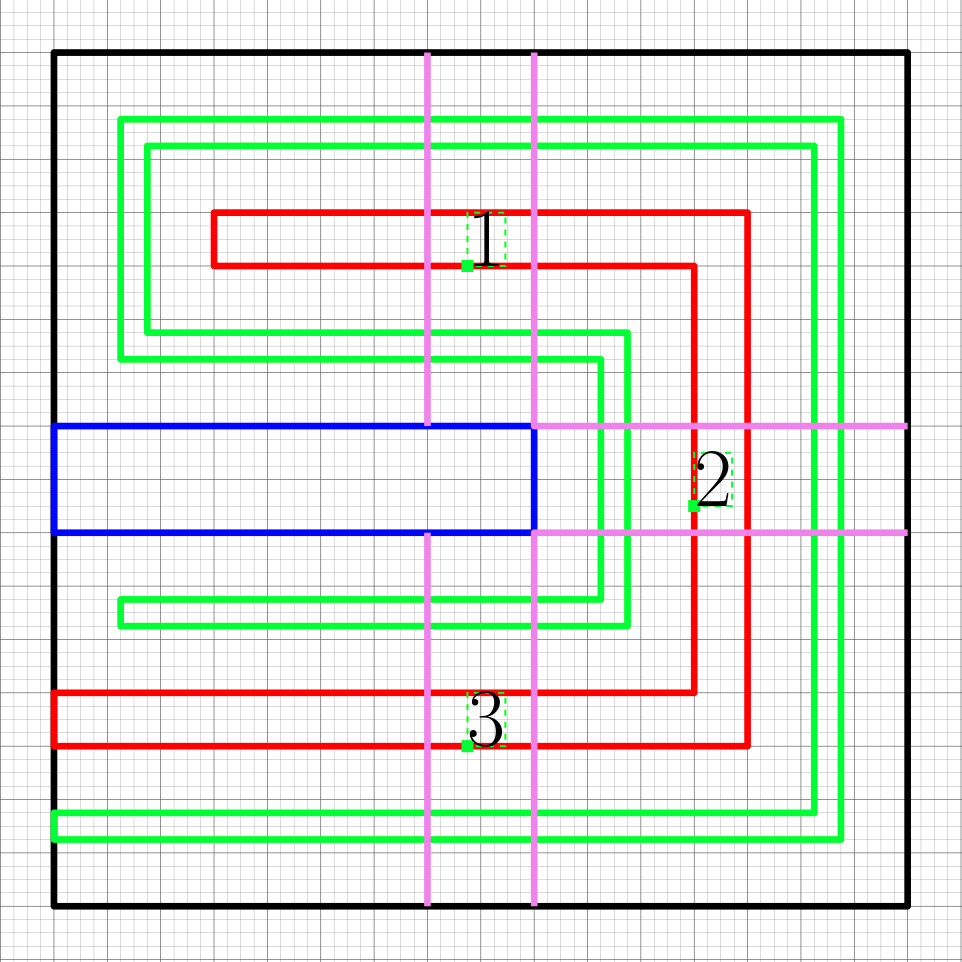}
              \caption{$\wa$ with $\mathcal{A}=16/7,7/3,12/5\dots$}
\end{figure}
      \end{minipage}
\end{minipage}
\medskip

The number of rectangles are bounded by the sum of the number of Turnings, of Terminals and of Separations, since there are at most one rectangle between two Turnings, between one Turning and one Terminal, etc.

\begin{Lem} Let $T(n)$ be the number of Turnings after the $n$-th digging. Then
\[
T(n) =\frac{9}{8}\cdot 2^n-4, \ \ \ n\geqslant 5
\]
\end{Lem}
\proof We calculate $T(n)$ by recurrence. Start with $T(1)=0, T(2)=2, T(3)=6, T(4)=14$. From Day $4$,
\begin{itemize}
\item each Straight produces no Turnings;
\item each Turning produces \underline{two Turnings};
\item one Terminal where we start our canal, produces one Separation and \underline{two Turnings};
\item the other Terminal produces \underline{two Turnings};
\item the Separation produces two Terminals and no Turnings.
\end{itemize}
Thus we get the recurrence formula
\[
T(n+1)=2T(n)+4, \ \ \ n\geqslant 4
\]
which implies that
\[
\aligned
\frac{T(n+1)}{2^{n+1}}&=\frac{T(n)}{2^n}+2^{1-n}\\
&=\frac{T(4)}{2^4}+(2^{1-n}+2^{2-n}+\dots+2^{1-4})\\
&=\frac{14}{2^4}+2^{-2}-2^{1-n},\ \ \ n\geqslant 4\\
T(n+1)&=\frac{9}{8}\cdot2^{n+1}-4, \ \ \ n\geqslant 4\\
T(n)&=\frac{9}{8}\cdot 2^n-4, \ \ \ n\geqslant 5
\endaligned
\]
\qed

Thus there are at most $T(n)+3=\frac{9}{8}\cdot 2^{n}-1$ rectangles on Day $n$ with $n\geqslant 5$. We can cover each rectangle by sub-squares with at most two overlaps.
\begin{figure}[H]
\begin{center}
              \includegraphics[width=0.2\linewidth]{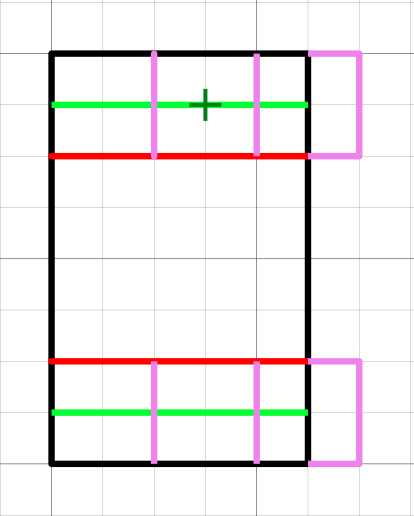}
              \end{center}
\end{figure}
So at each rectangle, the area that we covered twice are bounded by two sub-squares of length $t_{n}$. In total, we over-covered at most $2\big(T(n)+3\big)\,t_{n}^2$ area. Write $A(n)$ as the area of the island after $n$-th digging. Then
\begin{align}\label{eqn-area}
A(n)\leqslant \underbrace{N(t_{n})\,t_{n}^2}_{\text{total area covered by squares}}\leqslant A(n)+\underbrace{2(T(n)+3)t_{n}^2}_{\text{area of overlaps}}=A(n)+\big(\tfrac{9}{4}\cdot 2^n-\tfrac{9}{4}\big)\,t_{n}^2
\end{align}

\begin{Lem} The are of the island after $n$-th digging satisfies
\[
2^nt_n< A(n)<2^{n+1}t_n, \ \ \ n\geqslant 1
\]
\end{Lem}
\proof We calculate $A(n)$ by recurrence. Start with $A(0)=1$. After $n$-th digging,
\begin{figure}[H]
\begin{center}
              \includegraphics[width=0.4\linewidth]{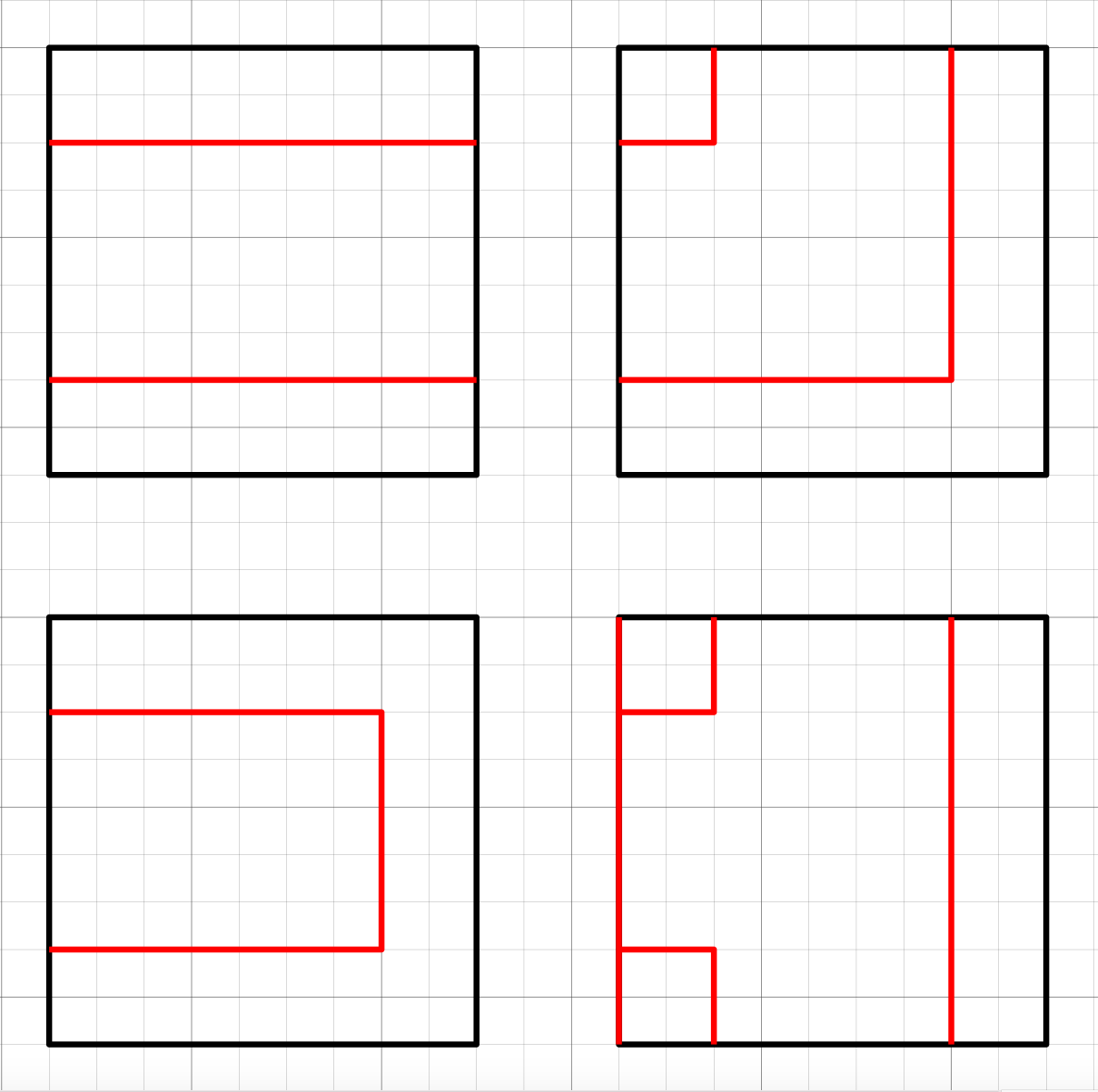}
              \caption{The change of the area from $A(n-1)$ to $A(n)$}
              \end{center}
\end{figure}
\begin{itemize}
\item the area of each Straight, no matter a square or a rectangle, is shrunk by $2/a_n$;
\item the area of each of the $T(n)$ Turnings is changed from $t_{n-1}^2$ to $t_{n-1}^2\cdot \frac{2}{a_n}$, i.e. shrunk by $2/a_n$;
\item the area of each of the two Terminals is changed from $t_{n-1}^2$ to $t_{n-1}^2\cdot\frac{1}{a_n^2}\cdot (3a_n-2)$;
\item the area of the Separation is changed from $t_{n-1}^2$ to $t_{n-1}^2\cdot\frac{1}{a_n^2}\cdot (a_n+2)$;
\end{itemize}
Thus
\[
\aligned
A(n)&=\left\{
\begin{aligned}
\Big(&A(n-1)-t_{n-1}^2\Big)\frac{2}{a_n}+t_{n-1}^2\cdot\frac{1}{a_n^2}\cdot (3a_n-2), \ \ \ n=1,2,3,\\
\Big(&A(n-1)-3t_{n-1}^2\Big)\frac{2}{a_n}+2t_{n-1}^2\cdot\frac{1}{a_n^2}\cdot (3a_n-2)+t_{n-1}^2\cdot\frac{1}{a_n^2}\cdot (a_n+2), \ \ \ n\geqslant 4
\end{aligned}
\right.\\
&=\frac{2}{a_n}A(n-1)+\frac{t_{n-1}^2}{a_n^2}(a_n-2), \ \ \ n\geqslant 1\\
&=\frac{2}{a_n}A(n-1)+t_n^2\,(a_n-2), \ \ \ n\geqslant 1\\
\frac{A(n)}{2^n\,t_n}&=\frac{A(n-1)}{2^{n-1}\,t_{n-1}}+\frac{t_n\,(a_n-2)}{2^n}, \ \ \ n\geqslant 1
\endaligned
\]
The last term is bounded by
\[
0< \frac{t_n\,(a_n-2)}{2^n} <\frac{t_n\,a_n}{2^n} =\frac{t_{n-1}}{2^n}
\]
Recall that $t_n=\prod\limits_{i=1}^na_i^{-1}<\prod\limits_{i=1}^n2^{-1}=2^{-n}$, so
\[
0< \frac{t_n\,(a_n-2)}{2^n} <2^{1-2n}.
\]
Thus for $n\geqslant 1$ 
\[
\aligned
\frac{A(n-1)}{2^{n-1}\,t_{n-1}}<\frac{A(n)}{2^n\,t_n}&<\frac{A(n-1)}{2^{n-1}\,t_{n-1}}+2^{1-2n},\\
A(0)<\frac{A(n)}{2^n\,t_n}&<A(0)+2^{-1}+2^{-3}+\dots\\
1<\frac{A(n)}{2^n\,t_n}&<1+1=2,\\
2^{n}t_n<A(n)&<2^{n+1}t_n.
\endaligned
\]
\qed

\begin{Def} Let $F(n),G(n)$ be two positive functions of $n\in\bz_{\geqslant 1}$. We write
\[
F(n)\simeq G(n),
\]
if there exist some constant $c>1$ such that
\[
c^{-1}\,G(n)<F(n)<c\,G(n).
\]
\end{Def}
 
So we can write $A(n)\simeq 2^n\,t_n$. In the equation~(\ref{eqn-area})
\[
\aligned
A(n)\leqslant N(t_{n})\,t_{n}^2&\leqslant A(n)+\big(\tfrac{9}{4}\cdot 2^n-\tfrac{9}{4}\big)\,t_{n}^2\\
&<2^{n+1}t_n+\tfrac{9}{4}\,2^n\,t_{n}^2\\
&<2^{n+1}t_n+\tfrac{9}{4}\,2^n\,t_{n}, \ \ \ (t_{n+1}<1)\\
&<(2+9/4)\,A(n).
\endaligned
\]
Thus $N(t_{n})t_{n}^2\simeq A(n)\simeq 2^n\,t_n$, implies $N(t_{n})\simeq \frac{2^n}{t_{n}}$ and
\[
\aligned
\dim_{\sf box}(\wa)&=\lim\limits_{n\rw+\infty}\frac{\ln\big(N(t_{n})\big)}{\ln(1/t_{n})}=\lim\limits_{n\rw+\infty}\frac{\ln(2^n/t_n)}{\ln(1/t_{n})}=\lim\limits_{n\rw+\infty}\frac{n\ln(2)-\ln(t_{n})}{-\ln(t_{n})}\\
&=1+\lim\limits_{n\rw+\infty}\frac{n\ln(2)}{-\ln(t_{n})}=1+\lim\limits_{n\rw+\infty}\frac{n\ln(2)}{\sum\limits_{i=1}^n\ln(a_i)}.
\endaligned
\]


\setlength\parindent{0em}
{\scriptsize Zhangchi Chen, Universit\'e Paris-Saclay, CNRS, Laboratoire de math\'ematiques d'Orsay, 91405, Orsay, France}\\
{\bf\scriptsize zhangchi.chen@universite-paris-saclay.fr}, {\bf\scriptsize https://www.imo.universite-paris-saclay.fr/$\sim$chen/}
\end{document}